\documentclass[11pt]{amsart}
\usepackage{amsmath,amssymb}
\usepackage{graphicx}

\usepackage{fullpage}
\usepackage{enumerate}
\usepackage{pgf,tikz}
\usetikzlibrary{arrows}
\usetikzlibrary{positioning}

\def\COMMENT#1{}
\let\COMMENT=\footnote
\def\TASK#1{}
\newcommand{\eps}{\varepsilon} 
\def\noproof{{\unskip\nobreak\hfill\penalty50\hskip2em\hbox{}\nobreak\hfill%
        $\square$\parfillskip=0pt\finalhyphendemerits=0\par}\goodbreak}
\def\endproof{\noproof\bigskip}
\newdimen\margin   
\def\textno#1&#2\par{%
    \margin=\hsize
    \advance\margin by -4\parindent
           \setbox1=\hbox{\sl#1}%
    \ifdim\wd1 < \margin
       $$\box1\eqno#2$$%
    \else
       \bigbreak
       \hbox to \hsize{\indent$\vcenter{\advance\hsize by -3\parindent
       \sl\noindent#1}\hfil#2$}%
       \bigbreak
    \fi}
\def\proof{\removelastskip\penalty55\medskip\noindent{\bf Proof. }}

\newtheorem{thm}{Theorem}[section]
\newtheorem{define}[thm]{Definition}

\newtheorem{lem}[thm]{Lemma}
\newtheorem{claim}[thm]{Claim}
\newtheorem{fact}[thm]{Fact}

\newtheorem{prop}[thm]{Proposition}

\newtheorem{question}[thm]{Question}

\newtheorem*{thm*}{Theorem}
\newtheorem*{define*}{Definition}
\newtheorem*{examp*}{Example}
\newtheorem*{lem*}{Lemma}
\newtheorem*{claim*}{Claim}
\newtheorem*{fact*}{Fact}
\newtheorem*{col*}{Corollary}
\newtheorem*{conj*}{Conjecture}

\begin{document}

\title{A discrepancy version of the Hajnal--Szemer\'edi theorem}

\author{J\'ozsef Balogh, B\'ela Csaba, Andr\'as Pluh\'ar and Andrew Treglown}
\thanks{JB: Department of Mathematics, University of Illinois at Urbana--Champaign, IL,  
USA, and Moscow Institute of Physics and Technology, Dolgoprodny,  Russian Federation, 
{\tt jobal@illinois.edu}. 
Research of this author is partially supported by  NSF Grants DMS-1500121, DMS-1764123, Arnold O. Beckman Research Award (UIUC) Campus Research Board 18132 and the Langan Scholar Fund (UIUC).
\\
\indent BC: Bolyai Institute, University of Szeged, Hungary, {\tt bcsaba@math.u-szeged.hu}. Research of this author was supported in part by the grant TUDFO/47138-1/2019-ITM of the
Ministry for Innovation and Technology, Hungary, and by NKFIH grant KH\_18 129597.\\
\indent AP: Department of Computer Science, University of Szeged, Hungary, {\tt pluhar@inf.u-szeged.hu}.\\
\indent AT: University of Birmingham, United Kingdom, {\tt a.c.treglown@bham.ac.uk}.}
\date{}
\begin{abstract}
A perfect $K_r$-tiling in a graph $G$ is a collection of vertex-disjoint copies of the clique $K_r$ in $G$ covering every vertex of $G$.
The famous Hajnal--Szemer\'edi theorem determines the minimum degree threshold for forcing a perfect $K_r$-tiling in a graph $G$.
The notion of discrepancy appears in many branches of mathematics. In the graph setting,
one assigns the edges of a graph $G$ 
 labels from $\{-1,1\}$, and one seeks substructures $F$ of $G$ that have `high' discrepancy (i.e. the sum of the labels of the edges in $F$ is far from $0$). In this paper we determine the minimum degree threshold for a graph to contain a perfect $K_r$-tiling of high discrepancy.
\end{abstract}
\maketitle

\section{Introduction} \label{Introduction}

\subsection{Discrepancy of graphs}
Classical discrepancy theory, or the study of irregularities of distribution, concerns with the following question: given some space, how evenly can one distribute a set of $n$ points in it (where here evenness is measured with respect to certain subsets)?
Perhaps the first significant result in the area is by Hermann Weyl on the criterion for a sequence to be uniformly distributed in
the unit interval. In the other direction, answering a question by van der Corput, van Aardenne-Ehrenfest proved that some irregularity of a point sequence in the unit interval is inevitable.
Since then discrepancy theory has become a widely studied area, with lots of ramifications and applications in ergodic theory, number theory, statistics, geometry, computer science, etc. For more details see the monograph by Beck and Chen~\cite{BeckChen}, the book by Chazelle~\cite{Ch} and the book chapter by Alexander and Beck~\cite{ABC}.

In this paper we study the discrepancy of graphs; a topic that lies in the wider framework of
\emph{hypergraph discrepancy theory} (see e.g.~\cite{BCsJP, Ch}). Before we can rigorously discuss this
topic we must introduce some definitions.
\begin{define}\label{defdisc}
Suppose $G$ is a graph and $f: E(G) \rightarrow \{-1,1\}$.
We say a subgraph $G'$ of a graph $G$ has \emph{discrepancy $t$ (with respect to $f$)} if $ \sum _{e \in E(G')} f(e) =t$ and \emph{absolute discrepancy $t$ (with respect to $f$)} if $\left | \sum _{e \in E(G')} f(e) \right | =t$.
\end{define}
If $G$ and $G'$ are $n$-vertex graphs, then we  say that $G$ contains a copy of $G'$ of \emph{high discrepancy}
(with respect to $f$) if there is a copy of $G'$ in $G$ with absolute discrepancy $\Omega (n)$.
Note that this concept also has a natural rephrasing in terms of Ramsey theory: given any $2$-colouring of the edges of $G$, one seeks a copy of $G'$ in $G$ whose edge set contains significantly more edges from one colour class than the other.

A natural question in graph discrepancy is to seek a fixed spanning subgraph $H$ of a graph $G$ 
of high discrepancy (or at least  discrepancy `far' away from zero). 
The first result of this type was obtained by Erd\H{o}s, F\"uredi, Loebl and S\'os~\cite{EFLS}:
they proved that, for some constant $c>0$, given any labelling $f: E(K_n) \rightarrow \{-1,1\}$ of $K_n$ and any 
 fixed spanning tree $T_n$ with maximum degree $\Delta$, $K_n$ contains a copy of $T_n$ of
absolute discrepancy at least $ c(n-1-\Delta)$. Note that in~\cite{EFLS} this result was phrased in the equivalent Ramsey setting.

In a previous paper~\cite{BCsJP}, Jing and the first three authors of this paper investigated the graph discrepancy problem of spanning trees, paths and Hamilton cycles
for various classes of graphs $G.$ For example, the following result determines the minimum degree threshold
for forcing a Hamilton cycle of high discrepancy.

\begin{thm}[Balogh, Csaba, Jing and Pluh\'ar~\cite{BCsJP}]\label{Oldthm:2.1}
Let  $0<c<1/4$  and $n\in \mathbb N$ be sufficiently large. If $G$ is an $n$-vertex graph with $$\delta(G)\geq(3/4+c)n$$
and $f: E(G) \rightarrow \{-1,1\}$, then 
there is a Hamilton cycle in $G$ with absolute discrepancy at least 
$ cn/32$ (with respect to $f$).
Moreover, if $4$ divides $n$, there is an $n$-vertex graph with $\delta(G)=3n/4$ and
an edge labelling $f: E(G) \rightarrow \{-1,1\}$ for which every Hamilton cycle has discrepancy $0$ (with respect to $f$).
\end{thm}
One can view such results about discrepancy 
as a measure of how robustly a graph contains
a spanning structure. Indeed, Theorem~\ref{Oldthm:2.1} implies that
every $n$-vertex graph $G$ with 
$\delta (G)> (3/4+o(1))n$ 
contains a Hamilton cycle that spans an
`unbalanced' collection of edges for any
partition $A \cup B$ of $E(G)$. (See~\cite{Sudakov} for a survey on other
measures of graph robustness.)

After submitting this paper, a multicolour extension of Theorem~\ref{Oldthm:2.1} was proven where the underlying graph is the random graph, see~\cite{newkri}.

\subsection{Perfect tilings in graphs}
An \emph{$H$-tiling} in a graph $G$ 
is a collection of vertex-disjoint copies of $H$ contained in $G$.
 An
$H$-tiling is \emph{perfect} if it covers all the vertices of $G$.
Perfect $H$-tilings are also often referred to as \emph{$H$-factors}, \emph{perfect $H$-packings} or \emph{perfect $H$-matchings}. 
$H$-tilings can be viewed as generalisations of both the notion of a matching (which corresponds to the case when $H$ is a single edge) and the Tur\'an problem (i.e. a copy of $H$ in $G$ is simply an $H$-tiling of size one).

Except for the case when $H$ contains no component of size at least $3$, the decision problem of whether a graph contains a perfect $H$-tiling is NP-complete (see~\cite{hell}). 
Thus, there has been substantial efforts to obtain sufficient conditions that force a graph to contain a perfect $H$-tiling.
In particular, 
a cornerstone result in extremal graph theory is the Hajnal--Szemer\'edi theorem~\cite{hs}, 
which characterises the minimum degree threshold that ensures a graph contains a perfect $K_r$-tiling. 

\begin{thm}[Hajnal and Szemer\'edi~\cite{hs}]\label{hs}
Every graph $G$ whose order $n$
is divisible by \(r\) and whose minimum degree satisfies $\delta (G) \geq (1-1/r)n$ contains a perfect $K_r$-tiling. Moreover, there are $n$-vertex graphs $G$
 with $\delta (G) = (1-1/r)n-1$ that do not contain a perfect $K_r$-tiling.
\end{thm}
There has also been much interest in the minimum degree threshold that ensures a perfect $H$-tiling for an arbitrary graph $H$.
After  earlier work on this topic (see~e.g.~\cite{alonyuster, kssAY}),  
K\"uhn and Osthus~\cite{kuhn, kuhn2}  determined, up to an additive constant, the minimum degree that forces a perfect $H$-tiling for \emph{any} fixed graph $H$.
Furthermore, there are now many different generalisations of the Hajnal--Szemer\'edi theorem. In particular, 
Kierstead and Kostochka~\cite{ore} proved an \emph{Ore-type} analogue,
Keevash and Mycroft~\cite{my} proved a version for \emph{$r$-partite graphs}, whilst there 
are now several generalisations of Theorem~\ref{hs} in the setting of \emph{directed graphs} (see e.g.~\cite{cdkm, forum}).

\subsection{Our main result}
In this paper we prove the following discrepancy version of the Hajnal--Szemer\'edi theorem.
\begin{thm}\label{main}
Suppose $r \geq 3$ is an integer and  let $\eta >0$. Then there exists $n_0 \in \mathbb N$ and $\gamma >0$ such that the following holds.
Let $G$ be a graph on $n \geq n_0$ vertices where $r$ divides $n$ and where
$$\delta (G) \geq \left (1-\frac{1}{r+1}+\eta \right )n.$$
Given any function $f: E(G) \rightarrow \{-1,1\}$ there exists a perfect $K_r$-tiling $\mathcal T$ in $G$ so that
$$\left | \sum _{e \in E(\mathcal T)} f(e) \right | \geq \gamma n.$$
\end{thm}
Comparing Theorem~\ref{main} with Theorem~\ref{hs} we see that having minimum degree just above that which forces a perfect $K_{r+1}$-tiling in fact ensures a perfect $K_r$-tiling of high discrepancy.
Moreover, the minimum degree condition in Theorem~\ref{main}  is essentially best-possible for all values of $r \geq 3$.
Interestingly, whilst the underlying extremal graph is the same for all $r\geq 3$ (the $(r+1)$-partite Tur\'an graph), the precise labelling of the edges we take is rather different depending on the value of $r$ modulo $4$. In Section~\ref{extremal} we construct  extremal labellings in the cases when $r \equiv 1,2 \ (\text{mod}  \ 4)$.
In the case when $r\equiv 0, 3 \ (\text{mod}  \ 4)$ the  extremal  labelling is easy to describe:
let $K$ be the complete graph $K_{r+1}$ with precisely half of its edges labelled with $1$, the remaining edges with $-1$ (the choice of $r$ ensures this is possible).
Then for any $n\in \mathbb N$ divisible by $r(r+1)$ consider the blow-up $G$ of $K$ with $n/(r+1)$ vertices in each class, and where the labellings of each edge in $G$ are induced by the labelling of $E(K)$.
It is easy to see that every perfect $K_r$-tiling in $G$ has discrepancy precisely $0$, whilst $\delta (G)=(1-1/(r+1))n$.
Moreover, in the case when $r(r+1)$ does not
divide $n$, the same construction $G$ is such that every perfect $K_r$-tiling has absolute discrepancy $O_r(1)$.

Note that the $r=2$ case (i.e.~perfect matchings) is covered by Theorem~\ref{Oldthm:2.1}. Indeed, it is easy to see that since the hypothesis of Theorem~\ref{Oldthm:2.1} forces a Hamilton cycle of high discrepancy, this ensures a perfect matching of high discrepancy. Moreover, consider the $4$-partite Tur\'an graph  $G$ on $n$ vertices
(where $4$ divides $n$). Label all edges incident to one of the vertex classes of $G$ with $-1$. All remaining edges are labelled $1$.
Then every perfect matching in $G$ has discrepancy $0$. Thus, perhaps surprisingly, this observation and Theorem~\ref{main} imply that the minimum degree threshold for forcing a perfect $K_3$-tiling of high discrepancy is the same as the analogous threshold for perfect matchings.
 
\smallskip

The paper is organised as follows. In the next section we introduce some notation and  definitions. In Section~\ref{extremal} we give the extremal examples for Theorem~\ref{main} in the cases when $r \equiv 1,2 \ (\text{mod}  \ 4)$. We introduce a number of tools that will be used in the proof of Theorem~\ref{main} in Section~\ref{useful}.
In Section~\ref{overview} we give an outline of the proof of Theorem~\ref{main} before giving the full proof in Section~\ref{proof}.
Finally, in Section~\ref{open} we present a number of open questions.

\section{Notation and definitions}
Let $G$ be a graph.
We write $V(G)$ for the vertex set of $G$, $E(G)$ for the edge set of $G$ and define
$|G|:=|V(G)|$ and $e(G):= |E(G)|$. Given a subset $X \subseteq V(G)$, we write $G[X]$ for the subgraph of $G$ induced by $X$ and 
$G\setminus X$ for the subgraph of $G$ induced by $V(G)\setminus X$.
The degree of $x $ is denoted by $d_G(x)$ and its neighbourhood
by $N_G(x)$. Given a vertex $x \in V(G)$ and a set $Y \subseteq V(G)$ we write $d _G (x,Y)$ to denote the number of edges $xy$ where $y \in Y$. 
Given a subgraph $F$  of $G$ we write
$d_G(x,F):=d_G(x,V(F))$.
Given disjoint vertex classes $X,Y\subseteq V(G)$,
we write $G[X,Y]$ for the bipartite graph with vertex classes $X$ and $Y$ whose edge set consists of all those edges in $G$
with one endpoint in $X$ and the other in $Y$;
 we write $e_G(X,Y)$ for the number of edges in $G[X,Y]$.

Suppose $G$ is a graph and $f: E(G) \rightarrow \{-1,1\}$.
 We say that $e \in E(G)$ is a \emph{$1$-edge} if $f(e)=1$ and a \emph{$(-1)$-edge} if $f(e)=-1$.
The \emph{$(-1)$-neighbourhood $N^-_G(x)$} of a vertex $x \in V(G)$ is the set of all vertices $y \in V(G)$
so that $xy $ is a $(-1)$-edge in $G$; the \emph{$1$-neighbourhood $N^+_G(x)$} of a vertex $x \in V(G)$ is the set of all vertices $y \in V(G)$
so that $xy$ is a $1$-edge in $G$.

The following notion of a $K_r$-template is crucial for the proof of Theorem~\ref{main}.
\begin{define}\label{template}
Let $F$ be a graph. A \emph{$K_r$-template of $F$ of size $s$} is a collection 
$\{H_1,\dots, H_s\}$ of not necessarily distinct copies of $K_r$ in $F$ for which there is some $s'\in \mathbb N$ so that every vertex $x \in V(F)$ lies in precisely $s'$ of the $H_i$.
(In fact, note we must have $s'=sr/|F|$.)
Suppose $f:E(F)\rightarrow \{-1,1\}$
and $\mathcal K:=\{H_1,\dots, H_s\}$ is a $K_r$-template of $F$.
We say that $\mathcal K$ has \emph{discrepancy $t$} if
$$\sum _{i=1} ^{s}\sum _{e\in E(H_i)} f(e)=t.$$
\end{define}
The following special labelled copies of $K_r$ appear in the proof of Theorem~\ref{main}.
\begin{define}\label{def1}
We write $K_r ^+$ for a copy of $K_r$ whose edges are each assigned $1$; define $K^- _r$ to be a copy of $K_r$ whose edges are each assigned $-1$.
The $(K_r,+)$-star is a copy of $K_r$ whose $1$-edges induce a copy of $K_{1,r-1}$. 
We call the root of this  $K_{1,r-1}$ the \emph{head} of the $(K_r,+)$-star.
We define the $(K_r,-)$-star and its head analogously.
\end{define}

We write $0<\alpha \ll \beta \ll \gamma$ to mean that we can choose the constants
$\alpha, \beta, \gamma$ from right to left. More
precisely, there are increasing functions $f$ and $g$ such that, given
$\gamma$, whenever we choose some $\beta \leq f(\gamma)$ and $\alpha \leq g(\beta)$, all
calculations needed in our proof are valid. 
Hierarchies of other lengths are defined in the obvious way.
Throughout the paper we omit floors and ceilings whenever this does not affect the
argument.

\section{The extremal examples}\label{extremal}
After its statement,  we introduced an extremal example for  Theorem~\ref{main} in the case when $r\equiv 0, 3 \ (\text{mod}  \ 4)$. In this section we first describe an extremal example for the case when $r \equiv 1 \ (\text{mod}  \ 4)$
and then give a construction for the $r \equiv 2 \ (\text{mod}  \ 4)$ case.
\begin{prop}\label{ex1} Let $m \in \mathbb N$, $r:=4m+1$ and $n \in \mathbb N$ be divisible by $2r(r+1)$.
Let $G$ be the complete balanced $(r+1)$-partite graph on $n$ vertices (and so $\delta (G)=(1-1/(r+1))n$). There is a function $f: E(G) \rightarrow \{-1,1\}$ so that for every perfect $K_r$-tiling $\mathcal T$ in $G$,
$\mathcal T$ has discrepancy zero (i.e.
$ \sum _{e \in E(\mathcal T)} f(e) =0$).
\end{prop}
\proof
Let $V_1,\dots, V_{r+1}$ denote the vertex classes of $G$; so $|V_i|=n/(r+1)$ for all $i \in [r+1]$.
Consider a copy $K$ of $K_r$ on vertex set $[r]$. Since $r=4m+1$ we can assign labels from $\{-1,1\}$ to each edge of $K$ so that the $(-1)$-edges induce a spanning $2m$-regular subgraph of $K$; the $1$-edges induce a spanning $2m$-regular subgraph of $K$.
Let $X,Y$ be a partition of $V_{r+1}$ so that $|X|=|Y|$.

We now define $f: E(G) \rightarrow \{-1,1\}$  as follows. The labelling of $K$ induces a labelling of the edges in $G':=G\setminus V_{r+1}$. That is, if $xy \in E(G)$ and $x \in V_i $, $y \in V_j$ where $1\leq i<j\leq r$, then $f(xy)=1$ if $ij$ is a $1$-edge in $K$;
 $f(xy)=-1$ if $ij$ is a $(-1)$-edge in $K$. Every vertex in $ X$ sends $1$-edges to each vertex in $V(G')$; every vertex in $Y$ sends $(-1)$-edges to each vertex in $V(G')$.

There are precisely three types of copy of $K_r$ in $G$:  \emph{Type~1} $K_r$ have every vertex  in $V(G')$; \emph{Type~2} $K_r$ have one vertex in $X$, the remaining vertices in $V(G')$; \emph{Type~3} $K_r$ have one vertex in $Y$, the remaining vertices in $V(G')$.
Note that a type 1 copy of $K_r$ has discrepancy $0$, a type 2 copy of  $K_r$ has discrepancy $r-1$, and a type 3 copy of $K_r$ has discrepancy $-r+1$.
Given any perfect $K_r$-tiling $\mathcal T$ in $G$, $\mathcal T$ must contain precisely the same number of type 2 and type 3 copies of $K_r$. Thus, $\mathcal T$ has discrepancy $0$, as desired.
\endproof
A similar function $f: E(G) \rightarrow \{-1,1\}$ to that in Proposition~\ref{ex1} now resolves the case when $r \equiv 2 \ (\text{mod}  \ 4)$. 
\begin{prop} Let $m \in \mathbb N$, $r:=4m+2$ and $n \in \mathbb N$ be divisible by $2r(r+1)$.
Let $G$ be the complete balanced $(r+1)$-partite graph on $n$ vertices (and so $\delta (G)=(1-1/(r+1))n$). There is a function $f: E(G) \rightarrow \{-1,1\}$ so that for every perfect $K_r$-tiling $\mathcal T$ in $G$,
$\mathcal T$ has discrepancy zero (i.e.
$ \sum _{e \in E(\mathcal T)} f(e)  =0$).
\end{prop}
\proof
Let $V_1,\dots, V_{r+1}$ denote the vertex classes of $G$; so $|V_i|=n/(r+1)$ for all $i \in [r+1]$.
Consider a copy $K$ of $K_r$ on vertex set $[r]$ whose edges are assigned labels from $\{-1,1\}$ so that there is precisely one more $1$-edge than $(-1)$-edge.
Let $X,Y$ be a partition of $V_{r+1}$ so that $|X|=\frac{(r-1)n}{2r(r+1)}$ and $|Y|=\frac{n}{2r}$.

We now define $f: E(G) \rightarrow \{-1,1\}$  as follows. As in Proposition~\ref{ex1}, the labelling of $K$ induces a labelling of the edges in $G':=G\setminus V_{r+1}$. Every vertex in $ X$ sends $1$-edges to each vertex in $V(G')$; every vertex in $Y$ sends $(-1)$-edges to each vertex in $V(G')$.

As before, there are precisely three types of copy of $K_r$ in $G$:  \emph{Type~1} $K_r$ have every vertex  in $V(G')$; \emph{Type~2} $K_r$ have one vertex in $X$, the remaining vertices in $V(G')$; \emph{Type~3} $K_r$ have one vertex in $Y$, the remaining vertices in $V(G')$.
Consider any perfect $K_r$-tiling $\mathcal T$ in $G$. Our aim is to show that  $\mathcal T$ has discrepancy $0$.

Note that  $\mathcal T$ contains precisely $\frac{n}{r}-\frac{n}{r+1}=\frac{n}{r(r+1)}$ copies of $K_r$ of type 1; each of these $K_r$s has discrepancy $1$.
Consider the subtiling $\mathcal T'$ of $\mathcal T$ induced by the type 2 and type 3 copies of $K_r$. Let $\mathcal T''$ be the $K_{r-1}$-tiling in $G'$ obtained from $\mathcal T'$ by removing all those vertices from $V_{r+1}=X\cup Y$.
Note that $\mathcal T''$ covers precisely $\frac{n}{r+1}-\frac{n}{r(r+1)}=\frac{(r-1)n}{r(r+1)}$ vertices in $V_i$ for each $i \in [r]$. In total $\mathcal T''$ consists of $n/(r+1)$ copies of $K_{r-1}$.
Moreover, for each pair $(i,j)$ with $1\leq i<j \leq r$, a precisely
$\frac{r-2}{r}$-proportion of the copies of $K_{r-1}$ in $\mathcal T''$ contain an edge $xy$ with $ x \in V_i$, $y \in V_j$.
Together, this implies that $\mathcal T''$ has discrepancy
$$\frac{r-2}{r} \times \frac{n}{r+1}.$$
Recalling that each edge incident to $X$ is a $1$-edge and each edge incident to $Y$ is a $(-1)$-edge, we conclude that $\mathcal T$ has discrepancy
$$\frac{n}{r(r+1)} + \left(\frac{r-2}{r} \times \frac{n}{r+1}\right ) + |X|(r-1) -|Y|(r-1)=\frac{n}{r(r+1)} +   \frac{(r-2)n}{r(r+1)}-\frac{(r-1)n}{r(r+1)}=0,$$
as required.
\endproof

\section{Useful results}\label{useful}
\subsection{The regularity lemma}\label{sec:reg}
In the proof of our main result we will use a discrepancy variant of Szemer\'edi's regularity
lemma~\cite{sze}. Before stating this result, we introduce some notation. 
The \emph{density} of a bipartite graph $G$ with vertex classes~$A$ and~$B$ is
defined to be
$$d(A,B):=\frac{e(A,B)}{|A||B|}.$$
Given any $\eps, d>0$, we say that $G$ is \emph{$(\eps,d)$-regular} if $d(A,B)\geq d$ and, for all sets
$X \subseteq A$ and $Y \subseteq B$ with $|X|\ge \eps |A|$ and
$|Y|\ge \eps |B|$, we have $|d(A,B)-d(X,Y)|< \eps$. 
We say that $G$ is \emph{$(\varepsilon, d)$-superregular} if all sets $X \subseteq A$ and $Y \subseteq B$ with $|X| \geq \varepsilon|A|$ 
and $|Y| \geq \varepsilon|B|$ satisfy that $d(X, Y ) > d$, that $d_G(a) > d|B|$ for all $a \in A$ and that $d_G(b) > d|A|$ for all $b \in B$.

Suppose $G$ is a graph with  edge labelling $f: E(G) \rightarrow \{-1,1\}$. Given disjoint $X,Y \subseteq V(G)$ we write
$G_{+}[X, Y]$  (or $(X,Y)^+_{G}$) for the bipartite graph with vertex classes $X,Y$ whose edge set consists of all those $1$-edges between $X$ and $Y$ in $G$.
We define $G_{-}[X, Y]$  and $(X,Y)^-_{G}$ analogously (now with respect to $(-1)$-edges).

We will apply the following variant  of Szemer\'edi's regularity lemma that can be easily deduced from the multicoloured version e.g. given in~\cite{ks}. 
\begin{lem} \label{reglem} For every $\varepsilon > 0$ and $\ell_0 \in \mathbb{N}$ there exists $L_0 = L_0(\varepsilon, \ell_0)$ so that the following holds.
Let $d \in [0, 1]$ and $G$ be a graph  on $n \geq L_0$ vertices with edge labelling $f: E(G) \rightarrow \{-1,1\}$. Then
there exists a partition $V_0, V_1, \ldots, V_\ell$ of $V(G)$ and a spanning subgraph $G'$ of $G$, such that the following conditions hold:
\begin{itemize}
\item [\rm (i)] $\ell_0 \leq \ell \leq L_0$;
\item [\rm (ii)] $d_{G'}(x) \geq d_G(x) - (2d + \varepsilon) n$ for every $x \in V(G)$;
\item [\rm (iii)] the subgraph $G'[V_i]$ is empty for all $1 \leq i \leq \ell$;
\item [\rm (iv)] $|V_0| \leq \varepsilon n$;
\item [\rm (v)] $|V_1| = |V_2| = \ldots = |V_\ell|$;
\item [\rm (vi)] for all $1 \leq i < j \leq \ell$  either $(V_i, V_j)^{+}_{G'}$ 
is an $(\varepsilon, d)$-regular pair or $G_{+}'[V_i, V_j]$ is empty;
\item [\rm (vii)] for all $1 \leq i < j \leq \ell$  either $(V_i, V_j)^{-}_{G'}$ 
is an $(\varepsilon, d)$-regular pair or $G_{-}'[V_i, V_j]$ is empty.
\end{itemize} \endproof
\end{lem}
We call $V_1, \dots, V_\ell$ {\it clusters}, $V_0$ the {\it exceptional set} and the
vertices in~$V_0$ {\it exceptional vertices}. We refer to~$G'$ as the {\it pure graph}.
The {\it reduced graph~$R$ of~$G$ with parameters $\varepsilon$, $d$ and~$\ell_0$} is the graph whose 
vertices are $V_1, \dots , V_\ell$ and in which $V_i V_j$ is an edge precisely when at least one of $(V_i,V_j)^+_{G'}$ and $(V_i,V_j)^-_{G'}$
is $(\varepsilon,d)$-regular. Associated with the reduced graph $R$ is an edge labelling $f_R: E(R) \rightarrow \{-1,1\}$ where
$f_R (V_iV_j):=1$ if $(V_i,V_j)^+_{G'}$ is $(\varepsilon,d)$-regular and $f_R (V_iV_j):=-1$ otherwise. (So if both $(V_i,V_j)^+_{G'}$ and $(V_i,V_j)^-_{G'}$
is $(\varepsilon,d)$-regular, then $f_R$ only `records' the former property.)

We will use the following well-known property of the reduced graph.
\begin{fact}\label{inherit}
Given a constant $c>0$, let $G$ be an $n$-vertex graph with $\delta (G)\geq c n$ 
that we   have applied Lemma~\ref{reglem} to  (with parameters $\varepsilon$, $d$ and~$\ell_0$).
Let R be the corresponding reduced graph. Then $\delta (R)\geq (c-2d-2\eps )|R|$. \endproof
\end{fact}
The following well-known result allows us to use subgraphs of a reduced graph as `templates' for
embedding in the original graph $G$.
\begin{lem}[Key lemma \cite{ks}]\label{keylem}
Suppose that $0 < \varepsilon < d$, that $q, t \in \mathbb{N}$ and that $R$ is a graph with $V(R) = \{v_1, \ldots, v_k\}$. 
We construct a graph $G$ as follows: Replace every vertex $v_i \in V(R)$ with a set $V_i$ of $q$ vertices and replace each edge of $R$ with an $(\varepsilon,d)$-regular pair. For each $v_i \in V(R)$, let $U_i$ denote the set of $t$ vertices in $R(t)$ corresponding to $v_i$.
Let $H$ be a subgraph of $R(t)$ with maximum degree $\Delta$ and set $h := |H|$. Set $\delta := d - \varepsilon$ and $\varepsilon_0 := \delta^{\Delta}/(2 + \Delta)$.
If $\varepsilon \leq \varepsilon_0$ and $t-1 \leq \varepsilon_0q$ then there are at least $$(\varepsilon_0 q)^h \ \mbox{labelled copies of $H$ in $G$} $$ so that if $x \in V(H)$ lies in $U_i$ in R(t), then $x$ is embedded into $V_i$ in $G$.
\end{lem}
The following fundamental result of Koml\'{o}s, S\'{a}rk\"{o}zy and Szemer\'{e}di~\cite{kssblowup}, known as the \emph{blow-up lemma}, essentially says that $(\varepsilon, d)$-superregular pairs behave like complete bipartite graphs with respect to containing bounded degree subgraphs.

\begin{lem}[Blow-up lemma \cite{kssblowup})] \label{blowup} Given a graph $F$ on vertices $\{1,\ldots, f\}$ and $d, \Delta > 0$, there exists an $\varepsilon_0 = \varepsilon_0(d,\Delta,f) > 0$ such that the following holds.
Given $L_1, \ldots , L_f \in \mathbb{N}$ and $\varepsilon \leq \varepsilon_0$, let $F^{*}$ be the graph obtained from $F$ by replacing each vertex $i \in F$ with a set $V_i$ of $L_i$ new vertices and joining all vertices in $V_i$ to all vertices in $V_j$ whenever $ij$ is an edge of $F$. 
Let $G$ be a spanning subgraph of $F^{*}$ such that for every edge $ij \in F$ the pair $(V_i,V_j )_G$ is $(\varepsilon, d)$-superregular. Then $G$ contains a copy of every subgraph $H$ of $F^{*}$ with $\Delta(H) \leq \Delta$.
\end{lem}

\subsection{An absorbing lemma}
We will apply the following well-known \emph{absorbing lemma} (which e.g. is a special case of~\cite[Theorem~4.1]{Tregs}).
Given a graph $G$ we say  a set $S \subseteq V(G)$ is a \emph{$K_r$-absorbing set for $Q \subseteq V(G)$}, if both
$G[S]$ and $G[S\cup Q]$ contain perfect $K_r$-tilings.
\begin{lem}\label{abslem}
Let $0 < 1/n \ll \nu \ll \eta  \ll 1/r$ where $n,r \in \mathbb N$ and $r \geq 2$. Suppose that $G$ is a graph on $n$ vertices with $\delta (G) \geq (1-1/r+\eta)n$.
Then $V(G)$ contains a set $M$ so that $|M|\leq \nu n$ and $M$ is a $K_r$-absorbing set for every $W \subseteq V(G) \setminus M$ such that $|W| \in r \mathbb N$ and  $|W|\leq \nu ^3 n$. 
\end{lem}

\section{Overview of the proof of Theorem~\ref{main}}\label{overview}
In the proof of Theorem~\ref{main} we will apply the regularity lemma to obtain the reduced graph $R$ of $G$
with an associated edge labelling $f_R: E(R) \rightarrow \{-1,1\}$.
Since the reduced graph $R$ `inherits' the minimum degree condition on $G$ (see Fact~\ref{inherit}), the Hajnal--Szemer\'edi theorem implies that
$R$ contains a perfect $K_{r+1}$-tiling $\mathcal T$.

In Claim~\ref{c1} we establish the following crucial property: (a) if  $\mathcal T$ has high absolute discrepancy (with respect to $f_R$), then we can use this structure in  $R$ as a framework to build a perfect $K_r$-tiling in $G$ with high absolute discrepancy (with respect to $f$). To build this tiling in $G$ we make use of the absorbing method.

We then establish another vital property of $R$: (b) if $R$ has a `small' subgraph $F$ so that $F$
has two $K_r$-templates with different discrepancies (with respect to $f_R$), then we can use this to again 
 build a perfect $K_r$-tiling in $G$ with high absolute discrepancy (see Claim~\ref{claim2}).
 
 We may therefore assume neither (a) nor (b) holds. This in turn forces the cliques of size at most $r+2$ in $R$ to have some very rigid structure. In particular, we deduce that every copy of $K_{r+1}$ in $R$ (therefore in our tiling $\mathcal T$) is one of the following: a  $K^+_{r+1}$; a $K^-_{r+1}$; a $(K_{r+1},+)$-star; a $(K_{r+1},-)$-star (see Claim~\ref{obs3}).

After this, we then argue that in fact almost all of the tiles in $\mathcal T$ are copies of $(K_{r+1},+)$-stars and $(K_{r+1},-)$-stars. Finally, we prove that there are two tiles $K,K'$ in $\mathcal T$
for which (b) must hold with respect to $F:=R[K\cup K']$, and so we do have 
 a perfect $K_r$-tiling in $G$ with high absolute discrepancy.

\section{Proof of Theorem~\ref{main}}\label{proof}
It suffices to prove the theorem in the case when $\eta \ll 1/r$.
Define additional constants $\gamma, \eps, d,\nu>0$ and $n_0, \ell_0, L_0 \in \mathbb N$ so that
\begin{align}\label{hierarchy}
0<1/n_0 \ll \gamma \ll 1/L_0 \leq  1/\ell _0  \ll \eps \ll d \ll \nu \ll \eta \ll 1/r.   
\end{align}
Here $L_0$ is the constant obtained
from Lemma~\ref{reglem} on input $\eps, \ell_0$.

Let $G$ be a graph on $n \geq n_0$ vertices as in the statement of the theorem.  Fix an arbitrary edge
labelling $f:E(G) \rightarrow 
\{-1,1\}$.

By Lemma~\ref{abslem} we obtain a  set of vertices $Abs\subseteq V(G)$ where $|Abs|\leq \nu n$ and where both $G[Abs]$ and $G[Abs \cup W]$ contain perfect $K_r$-tilings for any set $W \subseteq V(G)\setminus Abs$ of size at most
$\nu ^3n$ where $r$ divides $|W|$. Let $G_1:=G\setminus Abs$. 
Thus,
\begin{align}\label{min1}
    \delta (G_1) \geq \left (1-\frac{1}{r+1}+\frac{3\eta }{4} \right )n.
\end{align}


\subsection{Applying the regularity lemma}
Apply the regularity lemma (Lemma~\ref{reglem}) to $G_1$ with parameters $\eps, d, \ell_0$. 
We thus obtain clusters $V_1,\dots, V_{\ell}$ of size $m$ (where $\ell_0 \leq \ell \leq L_0$) an exceptional
set $V_0$ (of size at most $\eps n$) and a pure graph $G'_1$ of $G_1$. We may assume that $r+1$ divides $\ell$. (If not, we can achieve this by deleting at most $r$ of the clusters, and move the vertices in these clusters to the exceptional set $V_0$.) Further we obtain the reduced graph
$R$ of $G_1$ with an edge labelling 
$f_R: E(R) \rightarrow \{-1,1\}$ `inherited'
from $f$ (as defined in Section~\ref{sec:reg}). Note that (\ref{min1}) and 
Fact~\ref{inherit} imply that
\begin{align}\label{minr}
     \delta(R) \geq \left (1- \frac{1}{r+1} + \frac{\eta}{2} \right) \ell.
\end{align}
The following two claims will be used several times in our proof. The first implies that to obtain our desired perfect $K_r$-tiling in $G$
it suffices to find a perfect $K_{r+1}$-tiling in $R$ of high absolute discrepancy.
\begin{claim}\label{c1}
Suppose that $R$ contains a perfect $K_{r+1}$-tiling $\mathcal T_R$ with absolute
discrepancy $t\geq \eta ^2\ell$
(with respect to $f_R$). Then $G$ contains 
a perfect $K_r$-tiling with absolute discrepancy at least $\gamma n$ (with respect to $f$).
\end{claim}
\proof
Consider any copy $H$ of $K_{r+1}$ in $\mathcal T_R$.
Suppose that $H$ has discrepancy $t_H\in \mathbb Z$ (with respect to $f_R$).
The vertices $W_{1},\dots, W_{{r+1}}$ in $H$ are clusters in $G$. Write $G_H$ for the $(r+1)$-partite graph $G'_1[W_{1}\cup \dots \cup W_{r+1}]$.
Through repeated applications of the key lemma (Lemma~\ref{keylem}) we obtain that there  is a $K_r$-tiling
$\mathcal T_H$ in $G_H$ so that:
\begin{itemize}
    \item[(i)] All but precisely $\eps^{1/2} m$ vertices in $W_i$ are covered by $\mathcal T_H$ for each $i \in [r+1]$;
    \item[(ii)] Given any edge $xy \in E(\mathcal T_H)$, if $x \in W_i$ and $y \in W_j$ then
    $f(xy)=f_R(W_iW_j)$;
    \item[(iii)] Each copy of $K_r$ in 
    $\mathcal T_H$ 
    contains at most one vertex from every cluster $W_{i}$. Furthermore, given any $1\leq i<j\leq r+1$,
    a $\frac{r-1}{r+1}$-proportion
    of the $K_r$ in $\mathcal T_H$ contain an edge from $G'_1[W_i,W_j]$.
\end{itemize}
Note that (ii) follows from the definition of $f_R$; (iii) simply states that we embed copies of $K_r$ in $G_H$
in a balanced way, alternating which cluster $W_i$ is `uncovered by a copy of $K_r$'.
Since $H$ has discrepancy $t_H$, (ii) and (iii) imply that $\mathcal T_H$ has discrepancy
$$\frac{r-1}{r+1}\times |\mathcal T_H | \times t_H =\frac{(r-1)}{r}(1-\eps^{1/2})m t_H$$
(with respect to $f$).

Consider the $K_r$-tiling $\mathcal T'$ in $G'_1$ obtained by taking the union of the $\mathcal T_H$
for each $H$ in $\mathcal T_R$. 
By (i), $\mathcal T'$ contains all but  $|V_0|+\eps ^{1/2}m \ell \leq 2\eps ^{1/2} n$ of the vertices in $G_1$.
Noting that $\sum _{H \in \mathcal T_R } t_H\in \{t,-t\} $, we deduce that
$\mathcal T'$ has absolute discrepancy  
$$
\frac{(r-1)}{r}(1-\eps^{1/2})m t\geq \frac{2}{3}(1-\eps^{1/2}) \eta ^2 m\ell \geq \eta ^2 n/2
$$
(with respect to $f$).
Let $W$ be the set of vertices in $G_1$ uncovered by $\mathcal T'$; so $|W|\leq 2\eps^{1/2} n\leq \nu ^3 n$. Thus,
$G[Abs \cup W]$ has a perfect $K_r$-tiling $\mathcal T''$. As $|Abs \cup W|\leq 2 \nu n$,
$\mathcal T' \cup \mathcal T''$ is a perfect $K_r$-tiling in $G$ with absolute discrepancy at least
$\eta ^2n/2-\binom{r}{2}2 \nu n \geq \gamma n,$
as desired.
\endproof
The next claim gives us a useful condition that guarantees our desired perfect $K_r$-tiling in 
$G$; it will be used repeatedly through the proof.
\begin{claim}\label{claim2}
Let $F$ be a subgraph of $R$ on $p$ vertices where $r+1\leq p \leq 2r+2$.
Given some $s \leq r^{100}$, suppose that $F$ has two $K_r$-templates $\mathcal K=\{H_1,\dots, H_s\}$ and $\mathcal K'=\{H'_1,\dots, H'_s\}$, both
 of size $s$. If $\mathcal K$ and $\mathcal K'$
 have different discrepancies (with respect to $f_R$), then $G$ contains a perfect $K_r$-tiling
 with absolute discrepancy at least $\gamma n$.
\end{claim}
\proof
Let $W_1,\dots, W_p$ denote the clusters of
$G'_1$ that correspond to the vertices of $F$.
So if $W_iW_j \in E(F)$ and $f_R(W_iW_j)=1$
then $(W_i,W_j)^+ _{G'_1}$ is $(\eps,d)$-regular; otherwise if $W_iW_j \in E(F)$ and $f_R(W_iW_j)=-1$ then $(W_i,W_j)^- _{G'_1}$ is $(\eps,d)$-regular. A well-known property of regular pairs implies that we can delete $\eps ^{1/2}m$ vertices from each of these clusters to obtain
subclusters $W'_1,\dots, W'_p$ with the following properties:
 if $W_iW_j \in E(F)$ and $f_R(W_iW_j)=1$
then $(W'_i,W'_j)^+ _{G'_1}$ is $(2\eps,d/2)$-superregular;  if $W_iW_j \in E(F)$ and $f_R(W_iW_j)=-1$ then $(W'_i,W'_j)^- _{G'_1}$ is $(2\eps,d/2)$-superregular.
Write $m':=(1-\eps^{1/2})m$; so $|W_i'|=m'$ for
all $i \in [p]$.

Let $F^*$ be the  $p$-partite graph with
vertex classes $W'_1, \dots, W'_p$, and where for each $i \not = j$, there are all possible edges
between $W'_i$ and $W'_j$ precisely if $W_iW_j \in E(F)$; that is, $F^*$
is a blow-up of $F$. Define
$f_{F^*}:E(F^*)\rightarrow \{-1,1\}$ so that
$f_{F^*}(xy)=1$ if $x \in W'_i$, $y \in W'_j$
and $f_R (W_iW_j)=1$; $f_{F^*}(xy)=-1$ if $x \in W'_i$, $y \in W'_j$
and $f_R (W_iW_j)=-1$.

Write $t$ for the discrepancy of $\mathcal K$
and $t'$ for the discrepancy of $\mathcal K'$; by the assumption in the claim, $t \not=t'$. 
Note that we can use $\mathcal K$ as a `framework'
to find a perfect $K_r$-tiling $\mathcal T$ in $F^*$ as follows: consider any $H_k$ in $\mathcal K$ and let $W_{i_1},\dots, W_{i_r}$ be the vertices of $H_k$; in $\mathcal T$ there are $m'p/sr$ copies of $K_r$ corresponding to $H_k$ which contain precisely
one vertex from each of $W'_{i_1},\dots, W'_{i_r}$. Thus,  $\mathcal T $
has discrepancy  $m'pt/sr$ (with respect to 
$f_{F^*}$).

Similarly, we can use $\mathcal K'$ as a framework
to find a perfect $K_r$-tiling $\mathcal T'$ in $F^*$ of discrepancy  $m'pt'/sr$ (with respect to
$f_{F^*}$).

Now applying the blow-up lemma, this ensures
$G_0:=G'_1[W'_1\cup\dots \cup W'_f]$ contains 
two perfect $K_r$-tilings $\mathcal T_1$ and $\mathcal T_2$ with  discrepancy $m'pt/sr$ and  $m'pt'/sr$ respectively (with respect to
$f$). Note that 
$$|m'pt/sr -m'pt'/sr |\geq (1-\eps^{1/2})\frac{m}{s}
\geq
\frac{n}{2L_0r^{100}}\stackrel{(\ref{hierarchy})} {\geq} 2 \gamma n.  $$
Further, $G\setminus G_0$ comfortably satisfies
$$\delta (G\setminus G_0)\geq (1-1/r)n,$$
so contains a perfect $K_r$-tiling $\mathcal T_3$ by the Hajnal--Szemer\'edi theorem.
Therefore, both $\mathcal T_1 \cup \mathcal T_3$
and $\mathcal T_2 \cup \mathcal T_3$ are perfect $K_r$-tilings in $G$, whose discrepancies differ
by at least $2\gamma n$; thus, one of these 
perfect $K_r$-tilings has absolute discrepancy at least $\gamma n$, as desired.
\endproof
From now on we may assume that the hypotheses of  Claims~\ref{c1} and~\ref{claim2} fail; this will eventually lead to a contradiction, thereby proving the theorem.

\subsection{Properties of cliques in $R$}
The minimum degree condition on $R$ ensures the following easy observation.
\begin{claim}\label{obs1}
Let $1 \leq k \leq r+1$.
Every copy of $K_k$ in $R$ lies in a copy of $K_{r+2}$.
\end{claim}

We now use Claim~\ref{claim2} to prove that the copies of $K_{r+2}$ in $R$ have a limited number
of possible edge labellings.

\begin{claim}\label{obs2}
Every copy $K$ of $K_{r+2}$ in $R$ is one of the following: a $K^+_{r+2}$; a $K^-_{r+2}$; a $(K_{r+2},+)$-star; a $(K_{r+2},-)$-star.
\end{claim}
\proof
Consider an arbitrary Hamilton cycle $C$ in $K$. We obtain a  $K_r$-template $\mathcal K_C$ of $K$ of size $r+2$ by going around the Hamilton cycle as follows:
take each copy of $K_r$ whose vertices are $r$ consecutive vertices along $C$ and add it to $\mathcal K_C$. 

Consider any two Hamilton cycles $C=W_1\dots W_i W_{i+1} W_{i+2} W_{i+3} \dots W_{r+2}$ and $C'$ obtained from $C$ by reordering $W_i W_{i+1} W_{i+2} W_{i+3}$
as $W_i W_{i+2} W_{i+1} W_{i+3}$ (i.e. we just swap the order of $W_{i+1}$ and $W_{i+2}$).
Since we are assuming the hypothesis of Claim~\ref{claim2} does not hold, 
we must have that $\mathcal K_C$ and $\mathcal K_{C'}$ have the same discrepancy with respect to $f_R$.

This implies that $f_R (W_i W_{i+1})+f_R (W_{i+2} W_{i+3})=f_R (W_i W_{i+2})+f_R (W_{i+1} W_{i+3})$.
(The left hand side considers the contribution to the discrepancy of $\mathcal K_C$ not `present' in
the discrepancy of $\mathcal K_{C'}$; the right hand side considers the contribution to the discrepancy of $\mathcal K_{C'}$ not `present' in
the discrepancy of $\mathcal K_{C}$.)

The choice of the Hamilton cycle $C$ in $K$ was arbitrary. So this implies that 
\begin{align}\label{l1}
f_R (ab)+f_R (cd)=f_R (ac)+f_R (bd) \text{ for all distinct } a,b,c,d \in V(K).
\end{align}
Consider any $a \in V(K)$. Suppose $|N^-_K(a)|\geq 3$. Given any distinct $b,c,d \in N^-_K(a)$, (\ref{l1}) implies that $f_R(bd)=f_R (cd)$.
This implies that the edges in $N^-_K(a)$ are either all $1$-edges or all $(-1)$-edges. A similar argument holds if $|N^+_K(a)|\geq 3$.

In particular, this implies that if one of $N^-_K(a)$ and $N^+_K(a)$ is empty then $K$ is  one of the following: a $K^+_{r+2}$; a $K^-_{r+2}$; a $(K_{r+2},+)$-star; a $(K_{r+2},-)$-star.
We may therefore assume that both $N^-_K(a)$ and $N^+_K(a)$ are non-empty, and 
without loss of generality assume that  $|N^+_K(a)|\geq 2$.

Choose any distinct $c,d \in N^+_K(a)$ and $b \in N^-_K(a)$. Noting that $ac$ is a $1$-edge and $ab$ is a $(-1)$-edge, (\ref{l1}) implies $cd$ is a $1$-edge and $bd$ is a $(-1)$-edge. The choice of $c,d \in N^+_K(a)$ and $b \in N^-_K(a)$
was arbitrary so this implies all edges between $N^+_K(a)$ and $N^-_K(a)$ are $(-1)$-edges.

If $|N^-_K(a)|=1$ we are immediately done now: indeed, we have just argued that $b \in N^-_K(a)$ sends out $(-1)$-edges to everything else; and as $|N^+_K(a)|\geq 3$ in this case, all edges in $N^+_K(a)$ are $+1$-edges.
That is, $K$ is a $(K_{r+2},-)$-star.

Thus, we now may additionally assume $|N^-_K(a)|\geq 2$. Choose any distinct $c',d' \in N^-_K(a)$ and $b' \in N^+_K(a)$. Then (\ref{l1}) implies that $b'd'$ is a $1$-edge. This is a contradiction, as we already proved that all edges between  
$N^+_K(a)$ and $N^-_K(a)$ are $(-1)$-edges. Thus this case does not occur, and we are done.
\endproof

Combining Claims~\ref{obs1} and~\ref{obs2} we obtain the following.

\begin{claim}\label{obs3}
Let $1 \leq k \leq r+2$.
Every copy  of $K_{k}$ in $R$ is one of the following: a $K^+_{k}$; a $K^-_{k}$; a $(K_{k},+)$-star; a $(K_{k},-)$-star.
\end{claim}

\subsection{Using a perfect $K_{r+1}$-tiling in $R$}
Note that (\ref{minr}) and Theorem~\ref{hs} imply that $R$ contains a perfect $K_{r+1}$-tiling $\mathcal T$. By Claim~\ref{obs3}, there are only four types of $K_{r+1}$ in $\mathcal T$. Let $A$ denote the set of $K^+_{r+1}$ in $\mathcal T$;
 let $B$ denote the set of $K^-_{r+1}$ in $\mathcal T$; let $C$ denote the set of $(K_{r+1},+)$-stars in $\mathcal T$; let $D$ denote the set of $(K_{r+1},-)$-stars in $\mathcal T$.
Without loss of generality we may assume that 
\begin{align}\label{l2}
|B|+|C|\geq |A|+|D|.
\end{align}



\subsubsection{Assume that $A$ is non-empty.}
\begin{claim}\label{obs4}
Consider any vertex $V_a \in V(A)$ and any copy $K \in B$ of $K^-_{r+1}$. Then we may assume
$d_R(V_a,K) \leq r-2$ if $r$ is even; $d_R(V_a,K) \leq r-1$ if $r$ is odd.
\end{claim}
\proof
Write $K_A$ for the clique in $\mathcal T$ that contains $V_a$. Let $F:=R[K_A\cup K]$.

First consider the case when $r$ is even, and
suppose $V_a$ sends $r-1$ edges to $K$ in $R$.
Suppose $i$ of these edges are $1$-edges (and so $r-1-i$ of them are $(-1)$-edges). Let $X,Y \in V(K)$ be the vertices in $K$ that are not
incident to one of these $r-1$ edges.
We will prove that $F$ satisfies the hypothesis of Claim~\ref{claim2}.

Write $\mathcal K_A$ for the set of all copies
of $K_r$ in $K_A$, and $\mathcal K$ for the set of all copies
of $K_r$ in $K$;
so $|\mathcal K|=|\mathcal K_A|=r+1$.

Define $\mathcal K_1$ to be the $K_r$-template for $F$ of size $2r(r+1)$ that contains precisely
$r$ copies of each of the cliques in 
$\mathcal K_A \cup \mathcal K$. Note that
indeed $\mathcal K_1$ is a $K_r$-template for $F$
as each vertex $V \in V(F)$ is contained in
precisely $r^2$ of the cliques in $\mathcal K_1$. Since $K_A \in A$, $K \in B$, and 
$\mathcal K_1$ contains the same number of copies of cliques from $\mathcal K_A$ and $\mathcal K$,
$\mathcal K_1$ has discrepancy $0$ (with respect to $f_R$).

We define another $K_r$-template $\mathcal K_2$ for $F$ of size $2r(r+1)$ as follows:
\begin{itemize}
\item[(i)] for the clique $H \in \mathcal K_A$ that does not contain $V_A$,  add $2r-1$ copies of $H$ to 
$\mathcal K_2$;
\item[(ii)] for each clique $H \in \mathcal K_A$ that  contains $V_A$,
add $r-1$ copies of $H$ to 
$\mathcal K_2$;
\item[(iii)]  add to $\mathcal K_2$ $r$ copies of the clique in $F$ that contains $V_A$ and the $r-1$ vertices in
$V(K)\setminus \{X,Y\}$;
\item[(iv)]  add  $r+1$ copies of  each clique $H \in \mathcal K$ that contains both $X$ and $Y$;
\item[(v)] add one copy of each clique $H \in \mathcal K$ that  avoids one of $X$ and $Y$.
\end{itemize}
To prove that $\mathcal K_2$ is a $K_r$-template  for $F$ of size $2r(r+1)$ it suffices to prove that
every vertex $V \in V(F)$ lies in precisely $r^2$ of the cliques in $\mathcal K_2$:
if $V \in V(K_A)\setminus \{V_A \}$ then (i) and (ii) give that $V$ lies in
$(2r-1)+(r-1)(r-1)=r^2$ such cliques; (ii) and (iii) imply $V_A$ lies in 
$(r-1)r+r=r^2$ such cliques; if $V \in V(K)\setminus \{X,Y \}$ then (iii)--(v) imply that $V$ lies in
$r+(r+1)(r-2)+1\cdot 2=r^2$ such cliques; if $V \in  \{X,Y \}$ then (iv) and (iv) imply that 
$V$ lies in $(r+1)(r-1)+1\cdot 1=r^2$ such cliques.

To compute the discrepancy of $\mathcal K_2$ note that, compared to $\mathcal K_1$ it has:
one fewer clique from $\mathcal K_A$; $r-1$ fewer cliques from $\mathcal K$;
an additional $r$ cliques (from (iii)) that each have discrepancy $2i-\binom{r}{2}$.
As $\mathcal K_1$ has discrepancy $0$ this implies that $\mathcal K_2$ has discrepancy
$$-\binom{r}{2}+(r-1)\binom{r}{2}+r\left (2i -\binom{r}{2}\right )=2ir-r(r-1) \not =0$$
as $i \not =(r-1)/2$ (recall we assumed that $r$ is even).
So $F$ satisfies the hypothesis of Claim~\ref{claim2}.

Now suppose $r$ is odd and $V_a$ sends at least $r$ edges to $K$ in $R$.  We can fix $r-1$  such edges so that $i\not =(r-1)/2$ of them  are $1$-edges and $r-1-i$ of them are $(-1)$-edges.
Now arguing precisely as before we conclude $F$ satisfies the hypothesis of Claim~\ref{claim2} as desired.
\endproof

\begin{claim}\label{obs5}
Consider any $V_a \in V(A)$ and any  $K \in C$. Then  we may assume
$d_R(V_a,K) \leq r-2$ if $r$ is even; $d_R(V_a,K) \leq r-1$ if $r$ is odd.
\end{claim}
\proof 
Write $K_A$ for the clique in $\mathcal T$ that contains $V_a$. Let $F:=R[K_A\cup K]$.
Write $\mathcal K_A$ for the set of all copies
of $K_r$ in $K_A$, and $\mathcal K$ for the set of all copies
of $K_r$ in $K$;
so $|\mathcal K|=|\mathcal K_A|=r+1$.

The proof proceeds similarly to the previous claim. 
If $r$ is even, suppose $V_a$ sends $r-1$ edges to $K$ in $R$; if $r$ is odd
suppose $V_a$ sends $r$ edges to $K$ in $R$.
If all these edges avoid the
head\footnote{$K$ is a copy of a $(K_{k},+)$-star; the head of such a star was defined in Definition~\ref{def1}.}
$V_H$ of $K$ then we can argue precisely as in Claim~\ref{obs4}
to obtain two $K_r$-templates $\mathcal K_1$ and $\mathcal K_2$ of $F$, both with the same size, but different discrepancy. Note that how we construct $\mathcal K_1$ and $\mathcal K_2$  is identical to
the proof of Claim~\ref{obs4}, though the discrepancies will differ from that claim since now 
$K \in C$.

Next suppose $r$ is even and $V_a$ sends $r-1$ edges to $K$ in $R$, one of the endpoints being the head $V_H$.
Suppose $i$ of these edges are $1$-edges and $r-1-i$ of them are $(-1)$-edges.
Let $X,Y \in V(K)$ be the vertices in $K$ that are not endpoints of such edges.
Again, we choose  $\mathcal K_1$ and $\mathcal K_2$ as in Claim~\ref{obs4}.

That is, we define $\mathcal K_1$ to be the $K_r$-template for $F$ of size $2r(r+1)$ that contains precisely
$r$ copies of each of the cliques in 
$\mathcal K_A \cup \mathcal K$. We define
$\mathcal K_2$  as follows:
\begin{itemize}
\item[(i)] for the clique $H \in \mathcal K_A$ that does not contain $V_A$,  add $2r-1$ copies of $H$ to 
$\mathcal K_2$;
\item[(ii)] for each clique $H \in \mathcal K_A$ that  contains $V_A$,  add $r-1$ copies of $H$ to 
$\mathcal K_2$;
\item[(iii)]  add to $\mathcal K_2$ $r$ copies of the clique in $F$ that contains $V_A$ and the $r-1$ vertices in
$V(K)\setminus \{X,Y\}$;
\item[(iv)]  add  $r+1$ copies of  each clique $H \in \mathcal K$ that contains both $X$ and $Y$;
\item[(v)] add one copy of each clique $H \in \mathcal K$ that  avoids one of $X$ and $Y$.
\end{itemize}
The same argument as in Claim~\ref{obs4} implies both $\mathcal K_1 $ and $\mathcal K_2$ are 
$K_r$-templates for $F$ of size $2r(r+1)$.

To complete the proof we have to again show the discrepancies of $\mathcal K_1 $ and $\mathcal K_2$ are different.
Note that (i) and (ii) imply that $\mathcal K_2$ has  one fewer copy of $K^+_{r}$ from $\mathcal K_A$ compared to $\mathcal K_1$; compared to $\mathcal K_1$, $\mathcal K_2$ has an additional $r$ cliques arising from (iii);
from (iv) and (v) we conclude that $\mathcal K_2$ has $r$ fewer $(K_r,+)$-stars from $\mathcal K$ compared
to  $\mathcal K_1$; by (iv) $\mathcal K_2$ has one more copy of a $K^-_{r}$ from $\mathcal K$ compared
to  $\mathcal K_1$.
Thus, the difference in discrepancy between $\mathcal K_1$ and  $\mathcal K_2$ is precisely
$$
-\binom{r}{2}+r\left (2i+2(r-2)-\binom{r}{2}\right ) -r \left (-\binom{r}{2}+2(r-1)\right )
-\binom{r}{2}=2ri-r^2-r .
$$
As $r$ is even, this term is non-zero (since $i \not = (r+1)/2$ in this case).
Therefore, $\mathcal K_1$ and $\mathcal K_2$ are $K_r$-templates for $F$ of different discrepancies; that is,
the hypothesis of Claim~\ref{claim2} holds.

Next suppose $r\geq 5$ is odd, and $V_a$ has at least $r$ neighbours in $K$, including the head $V_H$.
 We can choose $r-1$  such neighbours, including $V_H$, so that $i$ of the corresponding edges incident to $V_a$ are $1$-edges  (and $r-1-i$ of them are $(-1)$-edges), where vitally, $i\not = (r+1)/2$. In particular, here we are using that
 $(r+1)/2<r-1$ to guarantee that we can choose $i$ as desired.
Then arguing as in the previous case we obtain two $K_r$-templates for $F$ of different discrepancies.
This argument also resolves the case when $r=3$ unless all the edges from $V_a$ to $K$ are $1$-edges; in which case we would be forced to `choose' $i=2=(r+1)/2$. However, in this case we have that
$V_a$ sends two $1$-edges to vertices in $V(K)\setminus \{V_H\}$. In this case can argue precisely as in Claim~\ref{obs4} to obtain 
two $K_r$-templates for $F$ of different discrepancies.
This completes the proof of the claim.



\endproof
By the last two claims we have that each $V_a \in V(A)$ has average degree of at most $r-1$ into each $K \in B\cup C$. Trivially $V_a$ has average degree of at most $r+1$ into each $K \in A \cup D$.
So by (\ref{l2}), each $V_a \in V(A)$ has average degree at most $r$ into each $K \in A \cup B\cup C \cup D$. This is a contradiction as $R$ has minimum degree $\delta (R) \geq (1-1/(r+1)+\eta /2)\ell$.
Thus we  conclude that $A$ is empty.

Further this implies  $B$ is small. Indeed, if $|B|\geq \eta ^2 \ell$  then  (\ref{l2}) implies that the perfect $K_{r+1}$-tiling $\mathcal T$ of $R$ has absolute discrepancy at least $\eta ^2 \ell$.
Thus the hypothesis of Claim~\ref{c1} holds, 
contradicting our assumption.

Therefore assume $A=\emptyset$ and $|B|\leq \eta ^2 \ell$. We now split into cases.

\subsubsection{Case 1: $r \geq 4$.}

Note that in this case we have
 $|D|-\eta^ 2 \ell  \leq |C|\leq |D|+\eta ^2 \ell$.
 Indeed, otherwise (\ref{l2}) 
 implies that the perfect $K_{r+1}$-tiling $\mathcal T$ of $R$ has absolute discrepancy at least $\eta ^2 \ell$.

Together with the fact that $\delta (R) \geq (1-1/(r+1)+\eta /2)\ell$, this implies that every $V_c \in V(C)$ has at least $(1-2/({r+1})+\eta /3)|C|$ neighbours in $D$. Thus, 
 this immediately implies
the following.
\begin{claim}\label{obs6}
 Given any $V_c \in V(C)$ there is some $K \in D$ such that $d_R(V_c, K) \geq r$. \qed
\end{claim}

Fix $V_c \in V(C)$ to be the head of some tile $K_C$ in $\mathcal T$.
So $V_c$ sends at least $r-1$ edges to $K\setminus \{V_H\}$ where $V_H$ is the head of $K$. Fix $r-1$ of these edges. Call the endpoints of these edges in $K$ \emph{good}.
Write $X$ for the vertex in $K\setminus \{V_H\}$ that is not good.
Write $\mathcal K_C$ for the set of all copies
of $K_r$ in $K_C$, and $\mathcal K$ for the set of all copies
of $K_r$ in $K$;
so $|\mathcal K|=|\mathcal K_C|=r+1$.

Set $F:=R[K_C \cup K]$.
Define $\mathcal K_1$ to be the $K_r$-template for $F$ of size $2r(r+1)$ that contains precisely
$r$ copies of each of the cliques in 
$\mathcal K_C \cup \mathcal K$. Note that
indeed $\mathcal K_1$ is a $K_r$-template for $F$
as each vertex $V \in V(F)$ is contained in
precisely $r^2$ of the cliques in $\mathcal K_1$. 

We define another $K_r$-template $\mathcal K_2$ for $F$ of size $2r(r+1)$ as follows:
\begin{itemize}
\item[(i)] for the clique $H \in \mathcal K_C$ that does not contain $V_c$,  add $2r-1$ copies of $H$ to 
$\mathcal K_2$;
\item[(ii)] for each clique $H \in \mathcal K_C$ that  contains $V_c$,
add $r-1$ copies of $H$ to 
$\mathcal K_2$;
\item[(iii)]  add to $\mathcal K_2$ $r$ copies of the clique in $F$ that contains $V_c$ and the good vertices;
\item[(iv)]  for each clique $H \in \mathcal K$ that  contains both $X$ and $V_H$,
add $r+1$ copies of $H$ to 
$\mathcal K_2$;
\item[(v)] add one copy of the clique $H \in \mathcal K$ that  avoids $X$;
\item[(vi)]  add one copy of the clique $H \in \mathcal K$ that  avoids  $V_H$.
\end{itemize}
To prove that $\mathcal K_2$ is a $K_r$-template  for $F$ of size $2r(r+1)$ it suffices to prove that
every vertex $V \in V(F)$ lies in precisely $r^2$ of the cliques in $\mathcal K_2$:
if $V \in V(K_C)\setminus \{V_c \}$ then (i) and (ii) give that $V$ lies in
$(2r-1)+(r-1)(r-1)=r^2$ such cliques; (ii) and (iii) imply $V_c$ lies in 
$(r-1)r+r=r^2$ such cliques; if $V \in V(K)\setminus \{X,V_H \}$ then (iii)--(vi) imply that $V$ lies in
$r+(r+1)(r-2)+1+1=r^2$ such cliques; if $V =X$ then (iv) and (vi) imply that 
$V$ lies in $(r+1)(r-1)+1=r^2$ such cliques; if $V =V_H$ then (iv) and (v) imply that 
$V$ lies in $(r+1)(r-1)+1=r^2$ such cliques.

We will now complete this case by showing that $\mathcal K_1$ and 
$\mathcal K_2$ have different discrepancies 
with respect to $f_R$; that is, the hypothesis of Claim~\ref{claim2} holds, as desired.

Write  $i$ for the number of  $(-1)$-edges in $F$
with one endpoint $V_c$, the other a good vertex.
So there are $r-1-i$ $1$-edges between $V_c$ and the good vertices.
Note that (i)  implies that $\mathcal K_2$ has  $r-1$ more copies of $K^-_{r}$ from $\mathcal K_C$ compared to $\mathcal K_1$; compared to $\mathcal K_1$, (ii) implies that $\mathcal K_2$ has $r$  fewer copies of  
$(K_r,+)$-stars from $\mathcal K_C$; the $r$ cliques from (iii) are contained in $\mathcal K_2$ but not
$\mathcal K_1$;
from (iv) and (v) we conclude that $\mathcal K_2$ has the same number of $(K_r,-)$-stars as $\mathcal K_1$;
by (vi) $\mathcal K_2$ has $r-1$ fewer copies of $K^+_{r}$ from $\mathcal K$ compared
to  $\mathcal K_1$.
Thus, the difference in discrepancy between $\mathcal K_1$ and  $\mathcal K_2$ is precisely
$$
-(r-1)\binom{r}{2}-r\left (-\binom{r}{2}+2(r-1)\right )
+
r\left (\binom{r}{2}-2i\right ) -(r-1) \binom{r}{2}=-r(r-1)-2ri<0.
$$
Therefore, $\mathcal K_1$ and $\mathcal K_2$ are $K_r$-templates for $F$ of different discrepancies; that is,
the hypothesis of Claim~\ref{claim2} holds, as required.



\subsubsection{Case 2: $r =3$.}
As $\delta (R) \geq (3/4+\eta /2)\ell$ and $|B|\leq \eta ^2 \ell$ we obtain the following.
\begin{claim}\label{obs7}
 Given any $V_c \in V(C)$ there is some $K \in C \cup D$ such that $d_R(V_c, K) =4$. \qed
\end{claim}
Fix $V_c \in V(C)$ to be the head of some tile $K_C$ in $\mathcal T$.
Write $\mathcal K_C$ for the set of all copies
of $K_3$ in $K_C$, and $\mathcal K$ for the set of all copies
of $K_3$ in $K$;
so $|\mathcal K|=|\mathcal K_C|=4$.
Set $F:=R[K_C \cup K]$.

{\it Subcase 2a: $K \in D$.}

Note that $V_c$ together with $K$ forms a copy of $K_5$ in $R$. As $K \in D$, Claim~\ref{obs2} tells us that the edge between $V_c$ and the head $V_H$ of $K$ is a $(-1)$-edge; all other edges between $V_c$ and $K$ are $1$-edges.

Define $\mathcal K_1$ to be the $K_3$-template for $F$ of size $24$ that contains precisely
$3$ copies of each of the cliques in 
$\mathcal K_C \cup \mathcal K$. Note that
indeed $\mathcal K_1$ is a $K_r$-template for $F$
as each vertex $V \in V(F)$ is contained in
precisely $9$ of the cliques in $\mathcal K_1$. 

We define another $K_3$-template $\mathcal K_2$ for $F$ of size $24$ as follows:
\begin{itemize}
\item[(i)] for the clique $H \in \mathcal K_C$ that does not contain $V_c$,  add $5$ copies of $H$ to 
$\mathcal K_2$;
\item[(ii)] for each clique $H \in \mathcal K_C$ that  contains $V_c$,
add $2$ copies of $H$ to 
$\mathcal K_2$;
\item[(iii)]  add to $\mathcal K_2$ one copy of each clique in $F$ that contains $V_c$ and precisely two of the vertices in $V(K)\setminus \{V_H\}$;
\item[(iv)]  for each clique $H \in \mathcal K$ that  contains  $V_H$,
add $3$ copies of $H$ to 
$\mathcal K_2$;
\item[(v)] add one copy of the clique $H \in \mathcal K$ that  avoids $V_H$.
\end{itemize}
It is easy to check that every $V \in V(F)$ lies in precisely $9$ cliques in $\mathcal K_2$; so
indeed $\mathcal K_2$ is $K_3$-template  for $F$ of size $24$.
Further, $\mathcal K_1$ has discrepancy $0$, $\mathcal K_2$ has discrepancy $-6$. Thus, the 
hypothesis of Claim~\ref{claim2} holds, as desired.



{\it Subcase 2b: $K \in C$.}

Note that $V_c$ together with $K$ forms a copy of $K_5$ in $R$. As $K \in C$, Claim~\ref{obs2} tells us that the edge between $V_c$ and the head $V_H$ of $K$ is a $1$-edge; all other edges between $V_c$ and $K$ are $(-1)$-edges.

We define $\mathcal K_1$ and $\mathcal K_2$ precisely as in Subcase 2a. That is,
define $\mathcal K_1$ to be the $K_3$-template for $F$ of size $24$ that contains precisely
$3$ copies of each of the cliques in 
$\mathcal K_C \cup \mathcal K$. Define  $\mathcal K_2$  as follows:
\begin{itemize}
\item[(i)] for the clique $H \in \mathcal K_C$ that does not contain $V_c$,  add $5$ copies of $H$ to 
$\mathcal K_2$;
\item[(ii)] for each clique $H \in \mathcal K_C$ that  contains $V_c$,
add $2$ copies of $H$ to 
$\mathcal K_2$;
\item[(iii)]  add to $\mathcal K_2$ one copy of each clique in $F$ that contains $V_c$ and precisely two of the vertices in $V(K)\setminus \{V_H\}$;
\item[(iv)]  for each clique $H \in \mathcal K$ that  contains  $V_H$,
add $3$ copies of $H$ to 
$\mathcal K_2$;
\item[(v)] add one copy of the clique $H \in \mathcal K$ that  avoids $V_H$.
\end{itemize}
In this subcase, $\mathcal K_1$ has discrepancy $0$, $\mathcal K_2$ has discrepancy $-12$. Thus, the 
hypothesis of Claim~\ref{claim2} holds, as desired.
This completes the proof of Theorem~\ref{main}.\qed

\section{Open problems}\label{open}
The \emph{$r$th power} of a Hamilton cycle $C$ is obtained from $C$ by adding an edge between every pair of vertices of distance at most $r$ on $C$.
The P\'osa--Seymour conjecture states that every $n$-vertex graph $G$ with minimum degree $\delta (G) \geq (1-1/(r+1))n$ contains the $r$th power of a Hamilton cycle.
Koml\'os, S\'ark\"ozy and Szemer\'edi~\cite{kss} proved this conjecture for sufficiently large $n$.

It is natural to seek a discrepancy analogue of the P\'osa--Seymour conjecture. We believe that the hypothesis of Theorem~\ref{main} additionally ensures
that the host graph $G$ contains the $(r-1)$th power of a Hamilton cycle with high discrepancy. 
Furthermore, the minimum degree in such a result should be best-possible (in the same sense Theorem~\ref{main} is best-possible). 
We believe the proof of such a result can be obtained via the connecting--absorbing method, and using Theorem~\ref{main} as a black-box (applied to the reduced graph of the host graph $G$); this would be a suitable project for a strong Master's student.
Note that such a result (combined with Theorem~\ref{Oldthm:2.1}) would show that $\delta (G)=(3/4+o(1))n$ is the threshold for a graph $G$ to contain \emph{both} a Hamilton
cycle of high discrepancy and the square of a Hamilton
cycle of high discrepancy.

It is also natural to seek an extension of Theorem~\ref{main} to perfect $H$-tilings for \emph{any} graph $H$.
\begin{question}
 Given any graph $H$, what is the minimum degree threshold that forces a perfect $H$-tiling  of high discrepancy in a graph $G$ (with respect to any edge labelling $f: E(G) \rightarrow \{-1,1\}$)?
\end{question}

A famous conjecture of Bollob\'as and Eldridge~\cite{bee}, and Catlin~\cite{cat}
asserts that every $n$-vertex graph $G$ with $\delta(G) \geq (rn-1)/(r+1)$ contains every $n$-vertex graph
$H$ with $\Delta (H)=r$. 
\begin{question}\label{ques2}
Given any $\eta>0$ and $r \geq 2$, does there exist an $n_0\in \mathbb N$ so that the following holds for all
$n \geq n_0$?
 Let $G$, $H$ be $n$-vertex graphs, and assume that $$\delta(G)\ge  (1-1/(r+2) + \eta)n,$$
where $r:=\Delta(H)$. Then  $G$ contains a copy of $H$ of high discrepancy  (with respect to any edge labelling $f: E(G) \rightarrow \{-1,1\}$).
\end{question}
Note that the Bollob\'as--Eldridge--Catlin conjecture has still not been fully resolved. So it seems extremely challenging to answer
Question~\ref{ques2} in general. 
However, our main result (Theorem~\ref{main}) resolves Question~\ref{ques2} in the affirmative when $H$ is a perfect $K_r$-tiling. It would be interesting to resolve Question~\ref{ques2} in cases for which the Bollob\'as--Eldridge--Catlin conjecture
is known to be true (in particular the case when $r \leq 3$).


\section{Acknowledgements}
The first and fourth authors are grateful to the BRIDGE strategic alliance between the University of Birmingham
and the University of Illinois at Urbana-Champaign. In particular, much of the research in this paper
was carried out whilst the fourth author was visiting UIUC.
The authors are also grateful to the referee for their helpful and careful review.

\end{document}